\documentclass{article}
\usepackage{amssymb, amsmath, latexsym}
\usepackage{epsfig}

\usepackage{psfrag}
\newtheorem{theorem}{Theorem}[section]
\newtheorem{definition}{Definition}[section]

\newtheorem{lemma}{Lemma}[section]
\newtheorem{claim}{Claim}[section]
\begin{document}
\large

\title{A knot bounding a grope of class $n$ is $\lceil \frac{n}{2}
\rceil$- trivial \footnote{Mathematics subject classification
57M25, 57M27}}
\author{James Conant}
\maketitle

\begin{abstract}

In this article it is proven that if a knot, $K$, bounds an imbedded grope
of class $n$, then the knot is $\lceil \frac{n}{2} \rceil$-trivial in the
sense of Gusarov and Stanford. That is, all type $\lceil \frac{n}{2}
\rceil$ invariants vanish on $K$. We also give a simple way to construct
all knots bounding a grope of a given class. It is further shown that this
result is optimal in the sense that for any
$n$ there exist gropes which are not
$\lceil
\frac{n}{2} \rceil+1$- trivial.
\end{abstract}

\section{Introduction}
\subsection{Origins}
Finite type invariants have been a hot topic of study in recent years,
having first been introduced in proto-form in a seminal paper of
Vassiliev[V], from which derives their alternative moniker ``Vassiliev
invariants''. Birman and Lin [BL], upon reading Vassiliev's paper were
able to give the by now familiar simple axiomatic condition for being a
finite type invariant of type n: Given a knot invariant $\nu$ taking
values in an abelian group extend it to knots with finitely many
transverse double points by the following formula, obligatory in any
paper on finite type invariants.
\begin{center}
\epsfig{file=obligatory,%
        height=1cm}
\end{center}
The invariant $\nu$ is \emph{finite type of type n} iff it vanishes on
knots with $n+1$ double points.  Birman and Lin also proved
that the coefficients of $x^n$ in the Jones polynomial under the change of
variables
$t\rightarrow e^x$ are type $n$ invariants. This is actually equivalent
to saying that the $n$th derivatives $J^{(n)}(1)$ are type $n$
invariants, and indeed this is used in the last section of the present
paper. 

\subsection{The work of Gusarov}
 Gusarov[G] takes a different tack, constructing a group of knots,
${\mathcal G}_n$, which is a quotient of the monoid of knots under
connected sum. The equivalence relation, as proven by Stanford and
Ng[NS], may be chosen to be that two knots are equivalent iff all
additive type n invariants are the same. An alternate description of this
group is given as follows. Given a knot
$K$ choose
$n+1$ disjoint groups of crossing changes $S=\{ s_1,\ldots, s_{n+1}\}$
for the knot. ($S$ is called a \emph{scheme} by Gusarov, or at least by
his translator.) If this scheme has the property that for some $L$,
$K_\sigma =L$ ($K_\sigma$ is the knot modified along the crossing changes
in $\sigma$.) for all nonempty $\sigma$, then we say $K\sim_n L$. ($K$ is
$n$-equivalent to $L$.)

Remarkably, $\sim_n$ is an equivalence relation [NS].
If we quotient the  monoid of knots under
$\#$ by $\sim_n$ we recover Gusarov's group
${\mathcal G}_n$. Denote elements in ${\mathcal G}_n$ by $[K]_n$, where
$K$ is a knot representing the equivalence class $[K]_n$. In fact,
for any scheme
$S$ the element
$Tot(K;S)\in {\mathcal G}_n$ is trivial where
$Tot(K;S):=\sum_{\sigma\subset S}(-1)^{|\sigma |} [K_\sigma]_n$, where
$|\sigma |$ is the cardinality of $\sigma$.  Indeed this is the main tool
of the present paper. This expands on the idea of Lin and Kalfagianni
[L-K] to just use the relation
$\sim_n$. Also, if we extend the above definition of finite type
invariants to links, this formula still holds in the following sense. Let
$\mu$ be a type $n$ invariant and $S$ a scheme of $n+1$ sets of crossing
changes of a link
$L$, then
\begin{gather}
\sum_{\sigma\subset S}(-1)^{|\sigma |}\mu (L_\sigma ) = 0. \label{altsum}
\end{gather}
\emph{[Proof]} (Following [G], Lemma 5.2)

An immediate consequence of the finite type axiom is the following:

If $S$ is a scheme of cardinality $n+1$ on $L$ where each $s_i\in S$ is a
single crossing change, then $\sum_{\sigma\subset S}(-1)^{|\sigma |} \mu
(L_\sigma )=0$, if $\mu$ is a type $n$ invariant. Our task is to prove
this when the $s_i$ contain more than one crossing change. We induct on
say $\sum_{i=1}^{n+1}|s_i|$. Given $S$ where $\sum |s_i| > n+1$, and
suppose without loss of generality that $s_1=s_1^\prime \cup
s_1^{\prime\prime}$ is a partition of $s_1$ into two nonempty sets. We
define 2 schemes of lower complexity:
\begin{gather*}
S^\prime =\{ s_1^\prime, s_2,\ldots ,s_{n+1}\} {\rm \hspace{.5em} on
\hspace{.5em}} L\\ S^{\prime\prime} = \{
\hat{s_1}^{\prime\prime},\hat{s_2},\ldots ,\hat{s}_{n+1} \}
{\rm\hspace{.5em} on\hspace{.5em}} L_{s_1^\prime}
\end{gather*}
where, if $s$ is a move on $L$, the move $\hat{s}$ denotes the induced
move on the link modified along $s_1^{\prime}$. 

\begin{gather*}
0+0=\sum_{\sigma\subset S^{\prime}}
(-1)^{|\sigma |}\mu (L_\sigma )
+\sum_{\sigma
\subset S^{\prime\prime}}(-1)^{|\sigma |}\mu ((L_{s_1^\prime})_\sigma )
=\\
=( \sum_{\tiny \begin{array}{c}\sigma\subset S^\prime \\
s_1^\prime \in \sigma\end{array}}(-1)^{|\sigma |}\mu (L_\sigma ) +
\sum_{\tiny \begin{array}{c}\sigma\subset S^{\prime\prime}\\
\hat{s_1^{\prime\prime}}\not\in\sigma\end{array}}(-1)^{|\sigma |}\mu
((L_{s_1^\prime})_\sigma )) + \\
\sum_{\tiny \begin{array}{c}\sigma\subset S^\prime \\
s_1^\prime \not\in \sigma\end{array}}(-1)^{|\sigma |}\mu (L_\sigma ) +
\sum_{\tiny \begin{array}{c}\sigma\subset S^{\prime\prime}\\
\hat{s_1^{\prime\prime}}\in\sigma\end{array}}(-1)^{|\sigma |}\mu
((L_{s_1^\prime})_\sigma )\\
\end{gather*}
\begin{gather*}
  =(\sum_{\tiny \begin{array}{c}\sigma\subset S^\prime \\
  s_1^\prime \in \sigma\end{array}}(-1)^{|\sigma |}\mu (L_\sigma ) +
  \sum_{\tiny \begin{array}{c}\sigma\subset S^{\prime}\\
  s_1^{\prime}\in\sigma\end{array}}(-1)^{|\sigma |+1}\mu
  ((L_{s_1^\prime})_\sigma ) ) +\\
  \sum_{\tiny \begin{array}{c}\sigma\subset S \\
  s_1 \not\in \sigma\end{array}}(-1)^{|\sigma |}\mu (L_\sigma ) +
  \sum_{\tiny \begin{array}{c}\sigma\subset S\\
  s_1\in\sigma\end{array}}(-1)^{|\sigma |}\mu
  ((L_{s_1^\prime})_\sigma )\\
=0 + \sum_{\sigma\subset S}(-1)^{|\sigma |}\mu (L_\sigma ) \Box
\end{gather*}

Note that this lemma generalizes the fact that $Tot( K;S)=0$ in two
senses: a) it holds for non-additive knot invariants and b) it holds for
\emph{link} invariannts.

 We'd like to point out also that a ``set of
crossing changes''
$s_i$ can be thought of a homotopy of the knot (or link) supported in a
disjoint union of balls. Indeed it is useful to think of it this way, in
which case a scheme $S$ is a set of ``disjointly supported homotopies.''
(Any homotopy of a knot beginning and ending with an embedding is
equivalent to a homotopy which is a set of disjointly supported finger
moves, i.e. crossing changes.)

\subsection{Gropes}
A grope, $G$, of class $n$, loosely, is a 2-complex representing an $n$
commutator [FT]. To define gropes recursively, however, we use a
different quantity, \emph{depth}. A depth
$1$ grope is defined to be a circle, while a depth 2 grope is defined as a
punctured surface. 
If you know what a depth
$<n$ grope is, to form a grope, $G$, of depth $n$, you take a
punctured surface and to each element of a prescribed symplectic basis you
glue a grope with depth
$<n$, such that at least one of these attached gropes is of depth $n-1$.

The class of a grope
$G$ is the term of the lower central series that the boundary circle
represents, or explicitly, if
$\{\alpha_i ,\beta_i\}$ is the symplectic basis and
$A_i, B_i$ are the gropes to be added, then ${\rm class}( G) =
\min_i\{{\rm class}(A_i)+{\rm class} (B_i)\}$. 

\subsection{Incorporating some geometry}
My result is then that if a knot bounds a grope of class $n$,
imbedded in
${\mathbb R}^3$, that that knot is trivial in ${\mathcal G}_{\lceil
\frac{n}{2}
\rceil}$. To do this I make repeated use of the fact that all the $Tot
(K;S)$'s are trivial in the group ${\mathcal G}_{|S|-1}$ by finding
appropriate collections of disjointly supported homotopies, the most
interesting of which come from the in/out trick 
defined in section 4. In a sense the main theorem is pretty easy to prove
if you don't mind not getting the optimal result. That is, without the
in/out trick, it is not so hard to prove that class $n$ gropes are
$\lfloor \frac{n}{2}\rfloor -1$- trivial. It is the in/out trick which
allows one to get those two extra groups of crossing changes for odd $n$.

It was originally suggested by Mike Freedman that class $n$ gropes might
always be $n-1$ trivial, (e.g. one group of crossing changes for every
`tip'.) This turns out to be overambitious by a factor of 2 and in the
last section we indeed deduce the existence of class $n$ gropes that are
not $\lceil \frac{n}{2}\rceil +1$-trivial.

An interesting consequence of the main theorem, (or the slightly weaker
one mentioned above) is that a knot bounding a grope of arbitrarily
large class cannot be distinguished from the unknot by finite type
invariants. It is a conjecture of Mike Freedman's that this phenomenon is
impossible. Indeed he conjectures that in any three manifold, you cannot
have an infinite imbedded grope, every stage of which is incompressible. 

\subsection{The work of X.S. Lin and E. Kalfagianni}
The main theorem of this paper is similar to that of a paper of X.S. Lin
and E. Kalfagianni[LK]. In that paper it is proven that knots which bound
certain immersed gropes of height $n+2$ are $l(n)$-trivial, where
$lim_{n\to \infty} l(n) = \infty$.  More specifically, they
consider immersed gropes such that all self-intersections occur away
from the bottom stage. There is also the restriction that the bottom
stage is \emph{regular}, which among other things implies that the
complement of the Seifert surface has free fundamental group. (It
should be noted that the obvious generalization of their and my result,
that all knots bounding immersed gropes are to some degree trivial,
meets with the problem that all knots bound immersed gropes of
arbitrary height, since the lower central series of a knot
complement stabilizes at
$\Gamma^2 =  {\mathbb Z}$.) 
The method of proof of their theorem, which precedes mine, is to find
crossing changes which implement the group-theoretic $n-1$-triviality of
an $n$-commutator, and as such is similar to the ideas presented in
the present paper, although the actual geometric implementation is
quite different. For instance their paper is mainly concerned with
the analysis of planar combinatorics, whereas mine uses more three
dimensional arguments.     

It turns out in their case, as well as mine, that one can not fully
realize the degree of triviality present in an $n$-commutator, the
problem being that one must be able to find $n$ geometric independent
moves which have the effect of deleting a letter in the commutator. For
the case of imbedded gropes one can show that only
$\lceil\frac{n}{2}\rceil+ 1$ of the moves are actually realizable
independently, whereas in the case of regular Seifert surfaces only
$\log_4 n/3 + 1$ of the moves are realizable independently.

\subsection{Synopsis}
In section 1 the introduction was given.

In section 2 we figure out how to put the grope $G$ into a nice form, and
using this form, to associate a decorated graph $\Gamma (G)$ to $G$.

In section 3 we reduce to the case when the bottom stage(The bottom stage is the one whose boundary is the knot
itself.)  is genus 1,
introducing two of the three types of moves (homotopies) we will need for
the main theorem.  We need to reduce to the genus 1 case in order to apply the
in/out trick which only works for gropes with bottom stage genus 1.

In section 4 we describe the in/out trick, and give some applications.
The trick is used in the proof of the main theorem and also in the
construction of the knots in section 6. 

In section 5, we finally polish off the main theorem. 

In section 6, we show that our result is optimal.

\vspace{1em}
{\bf Acknowledgements}

I wish to thank ARCS for their financial support. For the use of their
program, Knotscape, in calculating Jones polynomials I thank Morwen
Thistlethwaite and Jim Hoste. For many useful discussions and for the
invention of the in/out trick for seifert surfaces, I thank my advisor
Peter Teichner. For the original idea and for several stimulating
discussions I thank Mike Freedman. For help in sorting out a couple
minor details I thank Arthur
Bartels, Dave Bachman and Dave Letscher. For the original idea of showing
a knot is
$m$-trivial by looking for crossing changes, I thank X.S. Lin and Effie
Kalfagianni.

\newpage

\section{Standard position and the decorated graph $\Gamma (G)$}
\subsection{Standard position}
We begin by finding a nice handlebody surrounding the grope.  We need the
following definition of a particular 1 complex.

\begin{definition}
The 1-complex $\Xi_i$ is defined for all $i\in\mathbb{N}$ inductively
as follows. $\Xi_0$ is a point, while $\Xi_1$ is an interval. Now
suppose
$\Xi_{i-1}$ is defined and has 1-valent vertices $z_1,\ldots,z_k$. Form
$\Xi_i$ as the adjunction space gotten by gluing the midpoint of each of
$k$ intervals
$I_1,\ldots,I_k$ to the corresponding $z_1,\ldots,z_k$.
\end{definition} 
\begin{figure}
\begin{center}
\epsfig{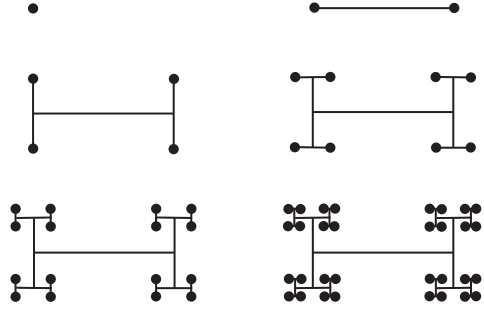}
\end{center}
\caption{The 1-complexes $\Xi_0$ to $\Xi_5$}
\end{figure}

Let $v_i$ be
imbedded circles representing the tips of the grope. For instance, if $G$
is a genus $2$ surface, there will be four $v_i$.

\begin{theorem}\label{standardposition}
  For every imbedded grope $G\subset \mathbb{R}^3$, there is a ball $B$
and handles $H_i\cong D^2\times I$ such that for all $t\in I$, the cross
section
$G\cap (D^2\times\{ t\}) \subset H_i$ is equal to $\Xi_{l(i)}$ for some
$l(i)$.  Further $v_i\cap (D^2\times\{ t\})$ is just a point in
$\Xi_1\subset
\Xi_{l(i)}$. Also $H_i\cap B = D^2 \times \partial I$. 
We also want there to be disks
$D_i\subset B$ where $\partial D_i = \gamma_i \cup \eta_i$ with
$\gamma_i \cap \eta_i = \partial \eta_i =\partial \gamma_i$ such
that
$\gamma_i\subset v_i$ and $\eta_i\subset \partial B$ and such that
$D_i\cap G=\gamma_i$. Also ${\rm int} D_i\cap {\rm int} D_j =\emptyset$ if
$i\neq j$. Finally, we require that
$B\cup_{i} H_i$ is a regular neighborhood of $G$. 
\end{theorem}

\emph{[Proof]}

First we show this is true for some model $ \mathcal{G}$ of $G$ in
$\mathbb{R}^3$.Once we do this, we are done. For if $f\colon\mathcal{G}\to
G$ then $f$ extends to give a PL-homeomorphism (or diffeomorphism
depending on which category you prefer) of regular neighborhoods
$\nu (\mathcal{G})\to \nu (G)$, which will transport the structure given
on the model.
A grope $G$ of depth 1 is just a circle. In this case we can let
$\mathcal G $ be an unknot. Take $B$ to be a small neighborhood of some
point of $\mathcal G$ and take the single handle $H_1$ to be a regular
neighborhood of the arc of $\mathcal G$ outside of $B$. The disk $D_1$ is
just a spanning disk of $\mathcal G$ intersected with $B$.

Now, for the inductive step suppose we have a grope $\mathcal G$ which is
formed as follows. Suppose the genus of the bottom stage of $\mathcal G$
is $g$ with a symplectic basis $\alpha_1,\ldots ,\alpha_g,\beta_1 ,\ldots,
\beta_g$ where $G$ is formed by attaching gropes $A_i, B_i$ to $\alpha_i,
\beta_i$ respectively. Since $A_i$ and $B_i$ have lower depth than
$\mathcal G$, the theorem is assumed to hold for them. So we have balls
$B_{A_i}, B_{B_i}$ together with handles $H_i(A_k), H_i(B_k)$ satisfying
the hypotheses of the theorem. Let $\delta_{A_k}$ be a small subarc of
$\partial A_k$ inside $B_{A_k}$ and $\epsilon_{A_k}$ a small arc joining
the endpoints of $\delta_{A_k}$ but with $\rm{int}(\epsilon_{A_k})$
contained in the interior of the bottom stage of $A_k$, such that
$\delta_{A_k}\cup\epsilon_{A_k}$ bound a disk $\Delta_{A_k}$. Modify
$B_{A_k}$ slightly, via a finger move disjoint from the rest of $\mathcal
G$ and the various $D_i$, so that $\partial (B_{A_k}\cup H_i(A_k))\cap A_k
=\epsilon_{A_k}$ with $\delta_{A_k}$ and $\Delta_{A_k}$ lying outside
$B_{A_k}$. Do this similarly for $B_k$. See figure \ref{gluean}.
\begin{figure}
\begin{center}
\epsfig{file=gluean,%
        height=4cm}
\end{center}
\caption{} \label{gluean}
\end{figure}

Now to form $\mathcal G$ attach annuli to the $A_k, B_k$ by gluing the
cores of the annuli to the boundaries of $A_k, B_k$ orthogonally to $A_k,
B_k$. Modify these by plumbing together the $A_k, B_k$ annuli for all $k$
and then connect summing all these together as in figure
\ref{glue2}, to form the genus $g$ bottom stage of the grope. Our new
handlebody
$B\cup_iH_i$ is formed as pictured in figure \ref{glue2}. The new handles
are the same as the old, but $B$ is formed by taking a small regular
neighborhood of $ G\backslash (\bigcup_i H_i \bigcup B_{A_k} \bigcup
B_{B_k})$.
\begin{figure}
\begin{center}
\epsfig{file=gluean2,%
        height=18cm}
\end{center}
\caption{}\label{glue2}
\end{figure}

If a handle looked like $\Xi_i\times I$ in $A_k$ then in $G$ it looks
like $\Xi_{i+1}\times I$, the effect of attaching an annulus. Also, the
$v_i$ for
$G$ are just made of all the $v_i$ for the $A_k$ and $B_k$ so they still
lie nicely in the handles as a subset of $\Xi_1\times I \subset
\Xi_{i+1}\times I$. As for the existence of disks $D_i$, consider a
handle $H$ for $G$ with tip
$v_H$. By hypothesis, there is a disk $D_H\subset B_{A_k}$ which extends
from
$v_H$ to
$\partial B_{A_k}$ and hence to $\partial B$. 
$\Box$   

This theorem gives a simple way of forming knots which bound gropes,
since we can imbed the handlebody in any way we please in ${\mathbb R}^3$.
Notice that the cores,
$v_i$, can be, using the disks
$D_i$, extended disjointly along annuli to curves $\overline{v}_i$ on 
$\partial (B\cup H_i)$.Now, proceeding with the advertised construction
of the graph
$\Gamma (G)$, we wish to group the cores $v_i$ into collections of cores
$V_i$,
$i=1,\ldots ,n$ where $n$ is the class of $G$. We want these $V_i$ to
have the property that if the collection of cores in some $V_i$ all bound
disks , $\Delta_{ij}$, into the complement of the grope, then the knot
$\partial G$ is isotopic to the unknot in a small regular neighborhood of
$G\cup_{i,j}\Delta_{ij}$. We do this inductively as follows. For a grope
with
$k(G)=2$, a Seifert surface, let
$V_1$ be formed by choosing one
$v_i$ from each pair of dual bands. $V_2$ is the set containing all the
other $v_i$. These obviously have the required property, since if $V_1$
bounds disks into the grope complement, surgery on these compressing
disks gives a spanning disk of $\partial G$. Now a grope with
$k(G)>2$, is formed by gluing gropes of lower depth, say $A_i, B_i$ to
a symplectic basis of the bottom stage,
$\alpha_i, \beta_i$. Suppose the class of $A_i$ is $s_i$ and that of
$B_i$ is $t_i$. Then by the inductive hypothesis we have groupings
$V_j^{\alpha_i}$, $j=1,\ldots,s_i$, and $V_j^{\beta_i}$, $j=1,\ldots
,t_i$. By definition the class of $G=\min_i \{ s_i+t_i\} =n$, say. Let
the elements of $\{ V_j^{\alpha_i}\}\cup \{ V_j^{\beta_i}\}$ be called
$\tilde{V}_1^i,\ldots ,\tilde{V}_{s_i+t_i}^i$. For $j<n$ define
$V_j^i=\tilde{V}_j^i$ and define $V_n^i=\cup
_{k=n}^{s_i+t_i}\tilde{V}_k^i$. Now define
$V_l :=
\cup_i V_l^i$. Now suppose
$V_l$ bounds disks into the grope complement. Then inductively for each
$i$, either
$A_i$ or $B_i$ can be surgered to produce a disk, since there exists a
$j$ such that
$V_j^{\alpha_i}\subset V_l^i$ or $V_j^{\beta_i}\subset V^i_l$
and hence for all $i$, $\partial A_i$ or $\partial B_i$ bounds
a disk. Hence a half basis of the bottom stage bound imbedded disks and so
surgery produces a spanning disk.

\begin{definition}
A set of handles has the \emph{trivializiation property} iff when caps of
these handles are abstractly added to the grope along the $v_i$ curves in
this set of handles, the grope becomes contractible. Another way to say
this is if the caps are added to the grope in a standard unknotted model
in
${\mathbb R}^3$, iterated surgery along the caps produces an unknotting
disk.
\end{definition}

 So now we
have a handlebody surrounding
$G$, with
$n$ groups of handles satisfying the trivialization property. Such a 
group
$V_i$ is said to be {\it framed unlinked} if the $\overline{v}_i$ bound
disks whose interiors intersect the grope only at handles not
associated to a core in $V_i$. This set of disks is called a
\emph{cap}. (When a disk does intersect a handle, by general position we
can assume it does so in a single level
$D^2\times
\{t\}
$.) If
$V_i$ is not framed unlinked, we say it is {\it framed linked}. The
reason for this terminology is that even if a group of handles $\{ H_i\}$
look like an unlink, a pushed out core $\overline{v}_i$ may link
with $v_i$ and hence will not be able to bound a disk into the grope
complement.

Fix a projection of the grope so that the 1-manifolds with
boundary,
$\overline{V}_i\cap B$ are standardly arranged in decreasing order as the
height function increases as in figure \ref{standheight}.
\begin{figure}
\begin{center}
\epsfig{file=standheight,%
        height=6cm}
\end{center}
\caption{}\label{standheight}
\end{figure}

To show this is possible, let $F: (\amalg I)\times I\rightarrow S^2$ be
an isotopy of $\cup(\overline{V}_i\cap B)$ to the standard picture
depicted in figure \ref{standheight}. Put a collar $C_1\cong S^2\times I$
on $B$ corresponding to the isotopy $F$. Let $C_2\cong S^2\times I$ be a
collar on
$B\cup C_1$ corresponding to a constant isotopy. Let
$C_3\cong S^2\times I$ be a collar on $B\cup C_1\cup C_2$ corresponding to
the isotopy inverse to
$F$.  We can think of the collar $C_1\cup C_2 \cup C_3$ as an
ambient isotopy of
$\cup \overline{V}_i {\hspace{.5em}\rm rel\hspace{.5em}} B$. We can now
take the new
$B$ to be
$B\cup C_1$, and the new handles to be regular neighborhoods of the
 part of the isotoped grope outside the new ball. The new disks
$D_i$ are formed by evolving the $\overline{V}_i$ outward along the radial
parameter.

\begin{definition}
A grope with the handlebody structure of theorem \ref{standardposition}
and designated groups of handles satisfying the trivialization property,
$V_i$, having a projection as described above, is said to be in
\emph{standard position.}
\end{definition}

\subsection{The decorated graph $\Gamma (G)$ associated to an imbedded
grope $G$.} \label{graphk}

\begin{definition}
Given a grope $G$ in standard position, we form a decorated graph $\Gamma
(G)$ as follows. We call the vertices $V_1,\ldots ,V_n$, corresponding to
the $n$ groups of handles satisfying the trivialization property. We put
an
$l$ next to a vertex if that group of handles is framed linked. We put an
edge between
$V_i$ and
$V_j$ ,$i< j$, if the group of handles $V_j$ ever cross over the group
$V_i$, with respect to the given projection.
\end{definition}  
 
As an exercise, note that if $\Gamma (G)$ consists of vertices with no
edges and no
$l$'s, then $\partial G$ is unknotted. This is because the cores all
bound disks, and also they are stacked with $V_1$ above $V_2$ above $V_3$,
etc. Thus, in particular, there is a plane separating $V_1$ and $V_2$
intersecting the ball in a level circle with respect to the height
function of the projection. So the disks bounding
$V_1$ say can be restricted to lie above
$V_2$ since if the disks ever ventured below the plane separating $V_1$
and
$V_2$ they could be surgered to lie above that plane, using the
3-manifold topologist's favorite tool, the inner-most disk argument.

\begin{definition}
Given a decorated graph $\Gamma$, the \emph{complexity}, $c(\Gamma )$ is
defined to be the number of edges, $E$, plus $\xi$, which is defined as
the number of vertices decorated with an $l$. That is $c(\Gamma )=E+\xi$.
\end{definition}

\begin{definition}
A group of vertices $V_{i_1}, \ldots , V_{i_k}$ is said to be
\emph{free}  if the $V_i$ are all framed unlinked and if for all $1\leq
s,t\leq k$ ($s\neq t$) there is no edge in the graph connecting $V_{i_s}$
with
$V_{i_t}$.
\end{definition}
\begin{figure}
\begin{center}
\epsfig{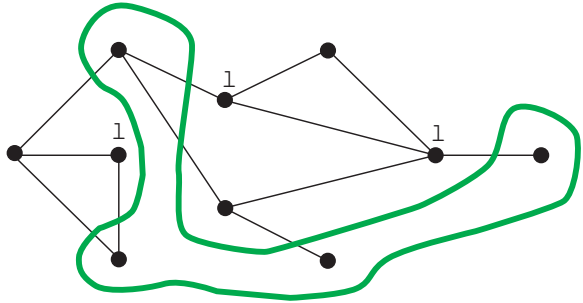}
\end{center}
\caption{The circled vertices are a group of 4 free
vertices.}\label{grapheg}
\end{figure}

\newpage
\section{Two types of moves and the reduction to the case where the
bottom stage of $G$ is genus 1}

\subsection{Move type I, a complexity reducing move}
Given an edge or an `l' in the graph $\Gamma$, we define a move which has
the effect of deleting the edge or `l', i.e. reducing $c(\Gamma )$.
Suppose the edge is between
$V_i$ and $V_j$. That means that some of the handles in $V_j$ cross over
or under some of the handles in $V_i$ in the wrong way. Then the move is
defined to be the homotopy which switches these handle crossings,
supported in balls associated to the crossings. (See [G] for instance.) In
order to remove an `l' from a vertex, suppose that vertex is $V_i$. To
unknot a handle in
$V_i$, first do handle crossings of the handle with itself so that the
handle bounds a disk which intersects only other handles. However we must
also make sure the handle is untwisted, which is to say that the pushed
out core
$\overline{v}_i$ of the handle bounds a disk which intersects only other
handles. Let the boundary of the disk that the handle bounds be the
longitude. Then Dehn twist to remove the appropriate number of multiples
of the meridian of the handle. This twist is supported in some small
 section of the handle
$D^2\times [a,b]$. Do this for every handle in $V_i$ to remove the `l'.

Notice that any number of type I moves may be performed simultaneously,
since the supports are by construction disjoint, with the effect that
the corresponding edges or `l's are deleted in $\Gamma$. 

\subsection{Move type II, moves on free vertices}

\begin{figure}
\begin{center}
\epsfig{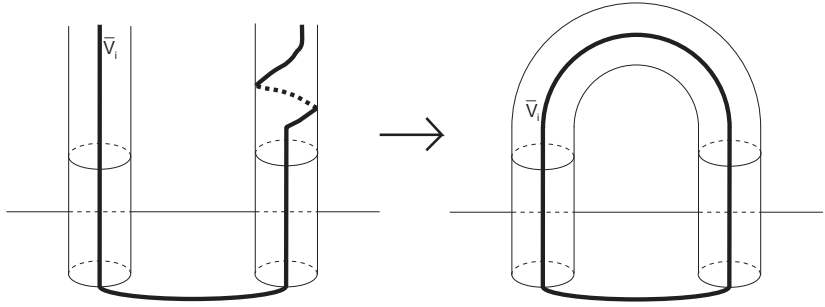}
\end{center}
\caption{A type II move trivializing a
handle.}\label{doublefivestar}
\end{figure}

Given a set, $F$, of $k$ free vertices we define $k$ moves as follows.
Since the vertices are free, there are planes in $S^3$ which separate the
groups of handles in $F$, and which intersect the ball of the grope's
handlebody in circles which lie standardly as level circles between the
attaching regions of the groups of handles. We can now choose homotopies
supported between the appropriate planes which contract the sets of
handles down to
  trivial handles within a small neigborhood of the ball as in figure
\ref{doublefivestar}. These moves obviously have disjoint support by
construction, and further doing any collection of them has the effect of
trivializing at least one set of handles with the trivialization
property. This has the effect of trivializing the grope. 

\subsection{The graph $\tilde{\Gamma}$.}
 Given a grope $G$, we
define a slightly different version of the graph defined in
section \ref{graphk}. Fix a projection of the handlebody where all the
$\overline{v}_i
\cap B$ occur in increasing order as height decreases. For the graph,
$\tilde{\Gamma}$ ,we let there be vertices $v_i$ for every handle in
the grope's associated handlebody, as opposed to one for each of the $n$
\emph{groups} of handles. We put an
$l$ next to the vertex if that core is framed linked in the previously
defined sense (since it is just one core you might say
\emph{framed knotted} instead), and we draw an edge between two vertices
if the corresponding handles cross in the wrong order in the projection. 

In terms of the graph $\tilde{\Gamma}$ we can still do type $\tilde{I}$
and type $\tilde{II}$ moves, defined in the obvious analogous way.
However, the result of doing a type $\tilde{II}$ move is no longer
neccessarily to trivialize the knot but instead to reduce the total genus
of the grope, where total genus is defined as the
sum of the genera of all the stages of the grope. (Since a trivialized 
handle has a core which bounds a disk, one
can iteratively surger along the successively produced disks as long as
the successive stages are genus one. When you hit a higher genus stage,
the surgery has the effect of lowering the genus of that stage by one.)
We have thus proved the following lemma.
\begin{lemma} \label{goodlemma}
If $\tilde{\Gamma}(G)$ has $k$ free vertices then $[\partial G]_{k-1}
= \sum \pm [\partial G_i]_{k-1}\in {\mathcal G}_{k-1}$, where $G_i$ is a
grope of lower total genus than $G$, but of the same class.
\end{lemma} 

\emph{[Proof]}

Let $S$ be the scheme of type $\tilde{II}$ moves defined above. Then
$\sum_{\sigma\subset S} (-1)^{|\sigma |} [\partial G_\sigma ]_{k-1}$. If
$\sigma\neq\emptyset$, then $G$ modified by $\sigma$ is of lower total
genus as analyzed above.
$\Box$.

\subsection{Genus 1 is sufficient}
Consider, toward a
contradiction, a counterexample which has minimal (total genus,
$c(\tilde{\Gamma})$), ordered lexicographically. This example has bottom
stage genus
$>1$, by assumption. Notice that 
$\tilde{\Gamma}$ has at least
$2n$ vertices, since for each pair of dual basis elements in the bottom
stage we get at least
$n$ vertices. I claim we can assume
$c(\tilde{\Gamma} )\leq \lceil \frac{n}{2}\rceil$. Otherwise, consider a
scheme,$S$, consisting of $\lceil \frac{n}{2} \rceil +1$ type $\tilde{I}$
moves.

By the triviality of the $Tot(K;S)$'s mentioned in the introduction,
inside the group ${\mathcal G}_{\lceil \frac{n}{2}\rceil}$ the knot
$K=\partial G$ is equivalent to a sum of knots of lower complexity and
equal total genus,
\begin{equation}
[K]_{\lceil \frac{n}{2}\rceil}=-\sum_{\emptyset\neq\sigma\subset S}
(-1)^{|\sigma |} [K_\sigma ]_{\lceil \frac{n}{2}\rceil}
\end{equation}

Each of these knots in the sum have reduced complexity, hence, by
minimality is $\lceil \frac{n}{2} \rceil$-trivial. Thus $[K]_{\lceil
\frac{n}{2}\rceil}=0$, contradicting that $K$ is a counterexample.

 So it suffices to consider knots with
$c(\tilde{\Gamma} )\leq \lceil \frac{n}{2}\rceil$. Now
\begin{eqnarray*}
\xi +E\leq \lceil \frac{n}{2}\rceil\\
\xi+2n-\lceil \frac{n}{2}\rceil\leq 2n-E\\
\Rightarrow \xi+\lceil \frac{n}{2}\rceil+1\leq \chi (\tilde{\Gamma}),
\end{eqnarray*}
 since $2n-\lceil \frac{n}{2}\rceil \geq \lceil \frac{n}{2} \rceil +1$.
On the other hand
$\chi (\tilde{\Gamma}) = b_0-b_1 =$
\# components -
\# cycles, implying there are at least $\chi$ components.
Hence there exist at least $\xi +\lceil \frac{n}{2}\rceil+1$ 
components, implying there are at least $\lceil \frac{n}{2}\rceil+1$
components without any framed linked vertices. We can choose a
free set of $\lceil \frac{n}{2}\rceil +1$ vertices by selecting
one vertex from each of these. So by
lemma \ref{goodlemma} , $[\partial G]_{\lceil\frac{n}{2}\rceil} = \sum \pm
[\partial G_i]_{\lceil \frac{n}{2}\rceil} = 0$ by the minimality
of total genus.
  But this is a contradiction, since
$K$ was supposed to be a counterexample.
We have thus proved the following lemma.

\begin{lemma}
If all class $n$ gropes with genus one bottom stage are $\lceil
\frac{n}{2}\rceil$-trivial, then {\emph all} class $n$ gropes are $\lceil
\frac{n}{2}\rceil$-trivial.
\end{lemma}
From now on assume the bottom stage of the grope is genus one.
\newpage
\section{Description of the In/Out Trick}
\subsection{Introduction}
  Whereas the two type of moves defined in the previous section
 preserve the grope structure, the move described
in this section, the in/out trick, does not. However, the move is
neccessary to prove the optimal result about grope triviality. Indeed we
also use it to construct the examples of section 6, for which lemma
\ref{4term}, section \ref{examples}, is needed. 

\subsection{The in/out trick} \label{random}
 Note that
$G$ divides naturally into two halves, the half attached to a particular
band of the bottom stage, and that attached to the dual band. Assume from
now on that the
$\overline{V}_i\cap B$ are in standard position such that all the
$\overline{V}_i\cap B$ on one half of the grope lie below all the
$\overline{V}_i\cap B$ on the other half.
In the handlebody, consider a framed unlinked $V_i=X$. Let
$\Delta_{x_1},\ldots \Delta_{x_m}$ be a cap. That is $\partial (\cup
\Delta_{x_i}) =
\overline{V}_i$, with the disks possibly intersecting the other handles.
If the cap does not intersect the other handles, then $\partial G$ is
unknotted.  

\begin{figure}
\begin{center}
\epsfig{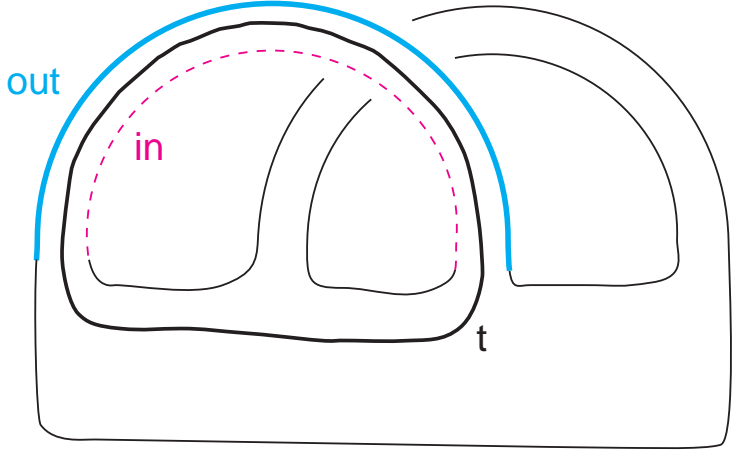}
\end{center}
\caption{The ``in'' and ``out'' arcs.}\label{inout1}
\end{figure}

 We consider two subarcs of
$K=\partial G$ called ``in'' and ``out'' by coloring the bottom stage of
$G$ as in figure
\ref{inout1}, where
$t$ is the curve to which
$X$'s half of the grope is attached. (In the inductive definition of the
$V_i$, it is obvious that $V_i$ lives in one half or the other.)

\begin{figure}
\begin{center}
\epsfig{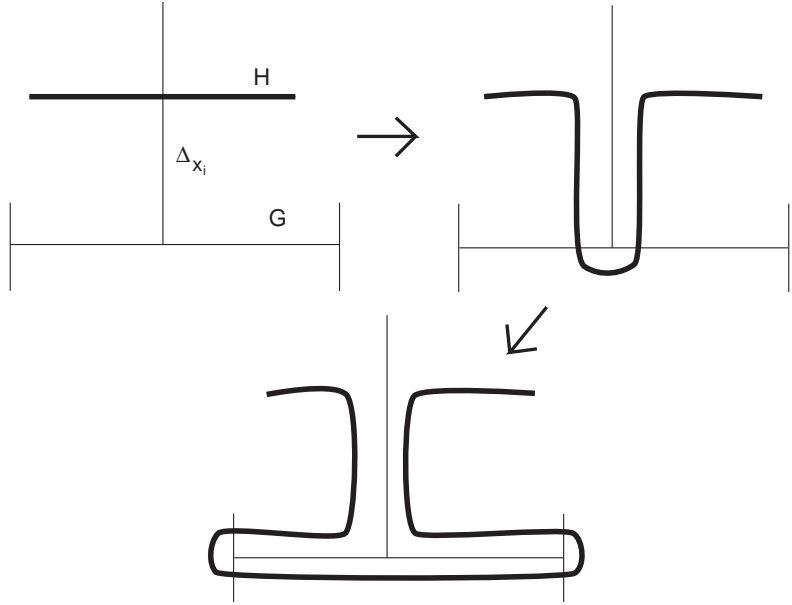}
\end{center}
\caption{Successive pushes of $H$. In the pictured case, the $X$ handle is
locally modelled on $\Xi_2 \times I$.}\label{inout2}
\end{figure}

Suppose $H$ is a handle intersecting $\Delta_{x_1}\cup\ldots\cup
\Delta_{x_m}$. Choose $m$ arcs inside the $\Delta_{X_i}$ from
$H\cap(\Delta_{x_1}\cup\ldots\cup\Delta_{x_m})$ to $X$. The endpoint of
each of these arcs lies on a handle of $V_i$ at some slice $D^2\times \{
t_0\}$. The grope slice at this point looks like some $\Xi_i$. Push
$H$ along these arcs as in figure \ref{inout2}. 

\begin{figure}
\begin{center}
\epsfig{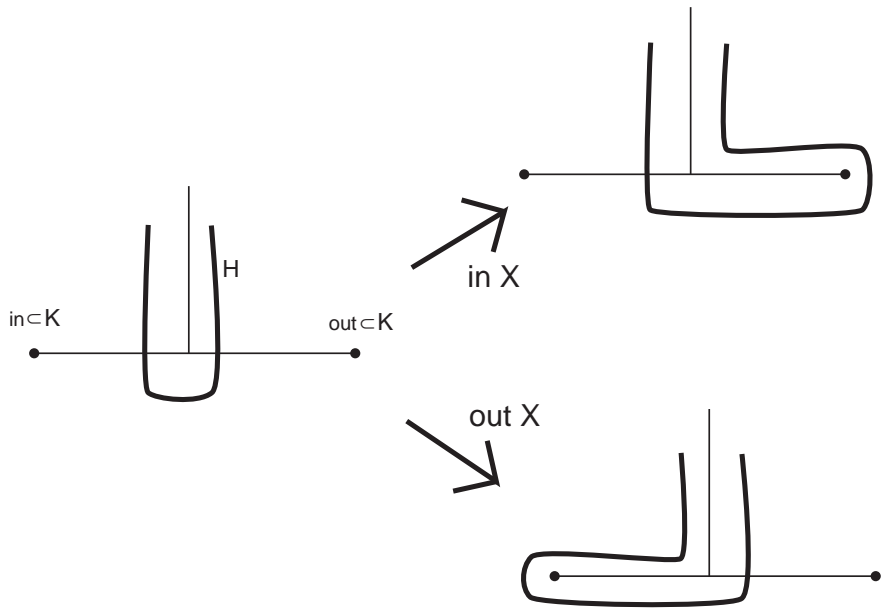}
\end{center}
\caption{The ``in'' and ``out'' moves.}\label{inout3}
\end{figure}

This introduces intersections of $H$ with a top stage of the grope.
(Although it's just an isotopy of $K$.) Continue pushing through
successive stages of the grope to eliminate the intersections, being sure
to push them down to the next stage in a small neighborhood of
$D^2\times
\{ t_0\}$. Continue doing this, for all handles, $H_i$, intersecting the
$\Delta_{x_i}$ until you've pushed to the bottom stage, but don't push
across the knot off the bottom stage (yet). 
If we push again, we'd be introducing actual crossing changes of the
knot. This preliminary isotopy will be called \emph{phase I} of the
in/out trick. We define the ``in'' and ``out'' moves as in figure
\ref{inout3}. These two homotopies are \emph{phase II}. They are clearly
disjointly supported after phase I.

Doing both ${\rm in} X$ and ${\rm out} X$ trivializes $K$, since it gives
a grope with
$X=V_i$ bounding disks. If we just do in $X$, then we can turn
the grope which $t$ bounds, $G^\prime$, into a disk, $\Delta$ in a regular
neighborhood of
$G^\prime\cup\cup_i \Delta_{x_i}$, by surgery. That is, glue in two
parallel copies of the $\Delta_{x_i}$ to make that stage of $G^\prime$ a
collection of disks. Iterate the procedure with these new disks until
$t$  bounds an (imbedded) disk, $\Delta$. (After all, we just removed all
intersections.) One subtlety is that this disk $\Delta$ will run
through the handles $H_i$, but this doesn't matter. Now, since we've
removed all intersections of $H$ between the ``in'' arc and $t$, we can
isotop the ``in'' arc along $\Delta$ to the arc $\mu$ as in figure
\ref{inout4}. This is \emph{phase III}.
\begin{figure}
\begin{center}
\epsfig{file=inout4,%
        height=7cm}
\end{center}
\caption{Doing ${\rm in}X$ gives the knot $t$ bounding the grope
$G^\prime$. }\label{inout4}
\end{figure}
 The
``out'' arc was never made to cross itself, so after the ``in'' arc
trivializes to $\mu$, the ``out'' arc can be isotoped back to its
original position. But now the band dual to $t$ pulls away, and we are
left with
$G^\prime$. This final isotopy is \emph{phase IV}. A similar
analysis holds for doing
${\rm out} X$, but one must pay attention to orientations. If $t$ is
oriented the same way as the ``in'' arc, then it will be oriented
oppositely to the ``out'' arc. Hence after doing ${\rm out} X$ we get the
knot $\rho (\partial G^{\prime} )$,where $\rho$ is the map reversing
orientation. (Note: it is not known whether finite type invariants can
ever distinguish a knot $K$ from $\rho(K)$.) 

For a genus one surface, the ``in'' and ``out'' arcs are symmetric so the
move $\rm{out} X$ gives the same (unoriented) result as $\rm{in} X$.
However, for a higher genus surface, the ``out'' move no longer works,
the problem occuring during phase IV, and this is why we need the bottom
stage of
$G$ to be genus one.

\subsection{Examples}\label{examples}
We now use the in/out trick to give a proof that every knot is 1-trivial.
This also follows from the main theorem and is well-known, but is good
for illustrative purposes.

\begin{figure}
\begin{center}
\epsfig{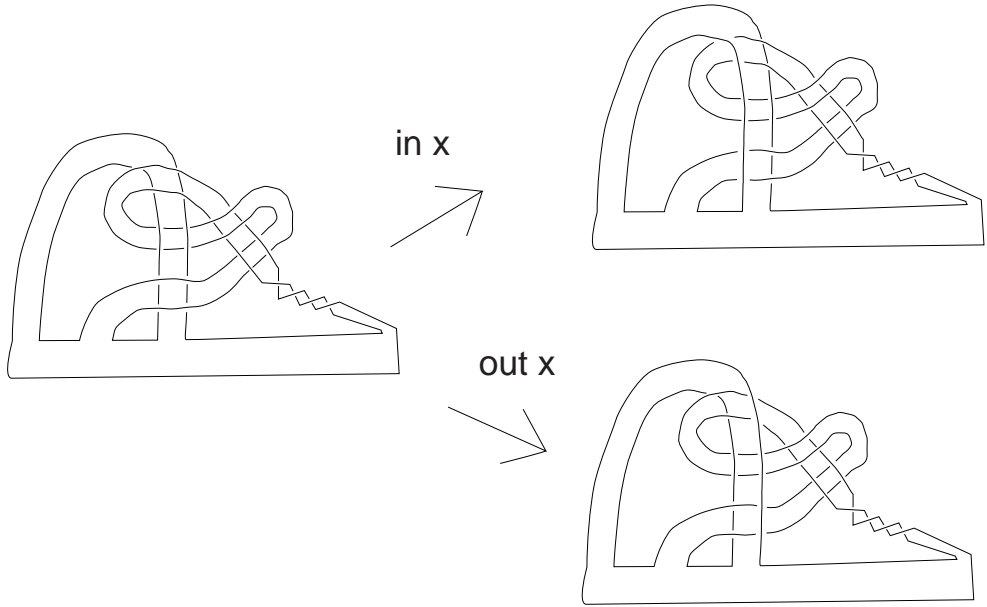}
\end{center}
\caption{The knots on the right are unknots.
}\label{inoutex}
\end{figure}

Suppose a knot, $K$, bounds a seifert surface with $k$ pairs of dual
bands $\{ x_i,y_i\}_{i=1}^k$. Consider the scheme $S=\{ s_1,s_2\}$ where
$s_1$ is the move which unknots and untwists the $x_1$ band and also does
crossing changes with other bands so that $x_1$ always crosses over them.
$s_2$ does a similar thing for $y_1$. Doing either $s_1$ or $s_2$ reduces
the genus of the seifert surface and we are done inductively. Doing both
gives a connected sum of a genus one knot that has unknotted bands and a
reduced genus knot. Thus it suffices to prove that a genus one knot with
unknotted bands, $x,y$, is $1$-trivial. But the scheme $\{ {\rm in} x,
{\rm out} x\}$ now trivializes the knot. In this simple case, ${\rm in} x$
(respectively ${\rm out}x$) may be visualized as the move making the
``in'' arc (respectively ``out'' arc) cross over everything in the
projection. See figure
\ref{inoutex}.

We conclude this section with an interesting calculation which will be
used in section \ref{sharp}.
\begin{lemma} \label{4term}
Consider a grope $G$ with genus one bottom stage which is formed by
gluing the gropes $G^\prime$ and $G^{\prime\prime}$ to the bottom stage.
They intersect in a point, $*$. 
\begin{figure}
\begin{center}
\epsfig{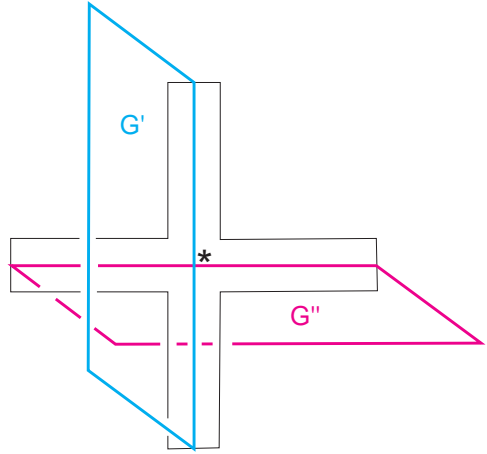}
\end{center}
\caption{The local picture at the bottom stage of $G$.
}\label{2inout}
\end{figure}
There are two ways to resolve this intersection inside the bottom stage
as pictured in figure \ref{3inout}.
\begin{figure}
\begin{center}
\epsfig{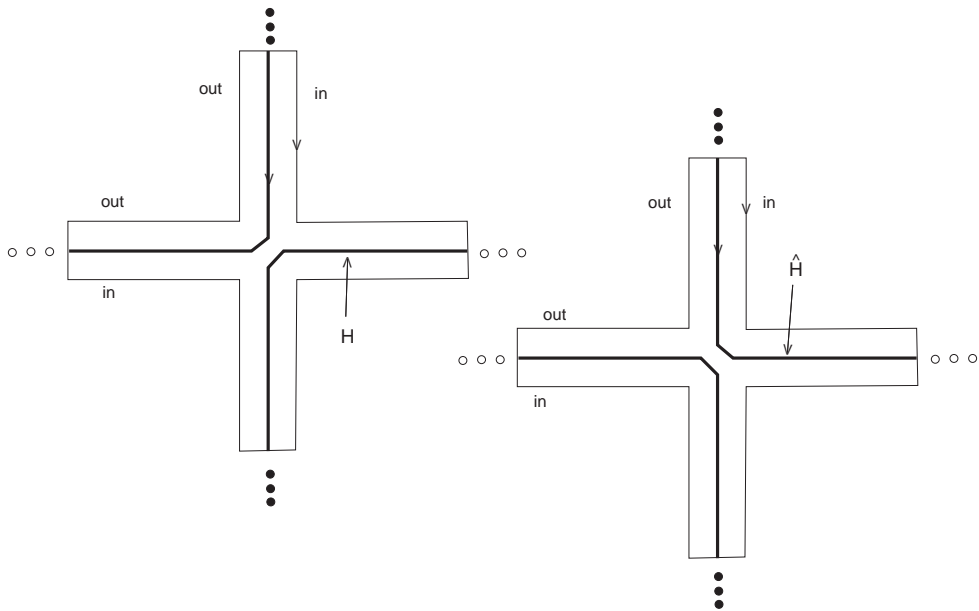}
\end{center}
\caption{The two resolutions.
}\label{3inout}
\end{figure}
These give rise to $2$ knots which are denoted $H$ and $\hat{H}$. Let $x$
be an unknotted vertex on the $G$ half and $y$ an framed unlinked vertex
on the
$G^{\prime\prime}$ half such that $\{ x,y\}$ is not an edge. Consider the
scheme $S =\{ {\rm in}x, {\rm out}x, {\rm in}y, {\rm out} y\}$. Then $Tot
(\partial G; S)=
 \sum_{\sigma \subset S} (-1)^{|\sigma |}\partial G_\sigma$, inside the
monoid ring ${\mathbb Z}Knots$, is equal to
$\partial G + H + \hat{H} + \rho (H) +
\rho (\hat{H} )$.
\end{lemma}
\emph{[Proof]}

\begin{figure}
\begin{center}
\epsfig{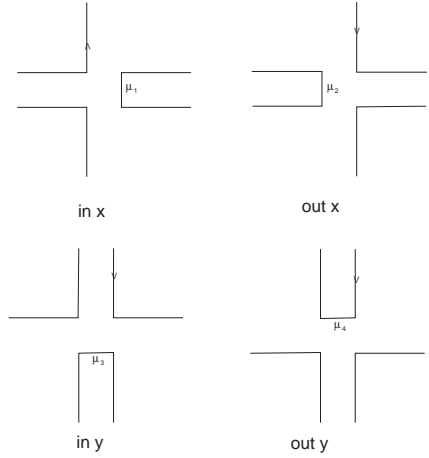}
\end{center}
\caption{Moves in S.
}\label{4inout}
\end{figure}

\begin{figure}
\begin{center}
\epsfig{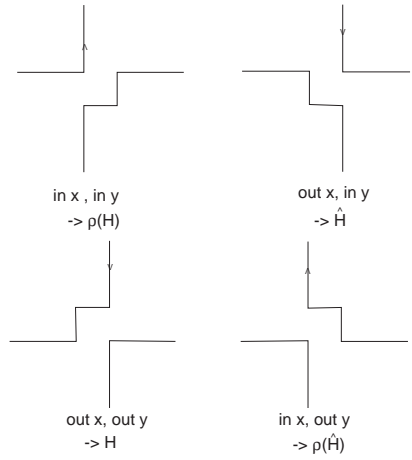}
\end{center}
\caption{Several pairs of moves in $S$.
}\label{5inout}
\end{figure}

\vspace{1em}
Consider figure \ref{2inout} depicting a neighborhood of $G^\prime\cap
G^{\prime\prime}$. Note the various moves in $S$ can be pictured as in
diagram
\ref{4inout}, the $\mu_i$ being the same as the $\mu$ arc previously
considered.
I claim the following:
\begin{gather*}
\sum_{\sigma\subset S}(-1)^{| \sigma |}\partial G_\sigma = (\partial G) -
(\partial G^\prime +\partial G^{\prime\prime}+\\
\rho(\partial G^\prime) +\rho (\partial G^{\prime\prime})) +(H+\rho (H)
+\hat{H}+\rho (\hat{H}))-(\partial G^\prime +\partial G^{\prime\prime} +
\rho (\partial G^\prime )+\rho (\partial G^{\prime\prime}))
\end{gather*}
which follows from the following facts: doing any single move in $S$ will
give the four terms of the second summand as was analyzed in section
\ref{random}. The third summand follows from diagram \ref{5inout} and the
fact that doing
${\rm in} x,{\rm out} x$ or ${\rm in} y,{\rm out} y$ together trivialize
the grope as analyzed in section \ref{random}. Of course some
justification is needed for diagram \ref{5inout}.

\begin{figure}
\begin{center}
\epsfig{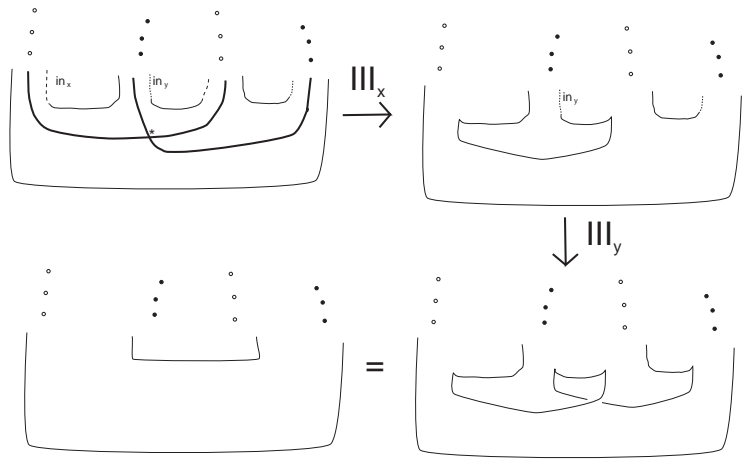}
\end{center}
\caption{Doing an in move on each half of the grope.}
\label{anotherinout}
\end{figure}

We must analyze what happens when we do, say, both ${\rm in}x$ and ${\rm
in}y$. Let $G^{\prime}_{I, II}$ (respectively $G^{\prime\prime}_{I, II}$)
be $G\prime$ (respectively $G^{\prime\prime}$) modified by phases I and
II of ${\rm in} x$ and ${\rm in} y$. Phase III of ${\rm in} x$
(respectively ${\rm in} y$) is supported in a regular neighborhood of
$G^\prime_{I,II} \cup x{\rm cap}$ (respectively $G^{\prime\prime}_{I, II}
\cup y{\rm cap}$). Note $(G^\prime_{I,II} \cup x{\rm cap}) \cap (
G^{\prime\prime}_{I, II}
\cup y{\rm cap})$ is the point $*$ in figure \ref{2inout}. Hence the phase
III isotopies are independent except near the end when the `in' arc gets
near $*$ soon to become the $\mu$ arc. So do the isotopies until they
are close to $*$ as in figure \ref{anotherinout}. But \ref{anotherinout}
is just a different picture of \ref{5inout}: ${\rm in}x, {\rm in}y$.

The fourth summand follows in the same way as the third, by considering
triplets of moves in $S$ and is left as an exercise to the reader.
Finally doing all moves in $S$ trivializes the grope.$\Box$

\section{The Main Theorem}

In this section, we prove the following
\begin{theorem} 
\label{mainthm}
Every class $n$ grope, $G$, is $\lceil \frac{n}{2} \rceil$-trivial.
\end{theorem}

\emph{[Proof]}

We may assume $n=2m+1$ since the even case follows by thinking of a class
$2m$ grope as a class $2m-1$ grope by forgetting a stage.
Also, we may assume $c(\Gamma )\leq m+1$ since we have $m+2$ moves in
hand to reduce complexity. Now a set of $m$ free vertices exists by the
following euler characteristic argument. ($b_i$ denote Betti numbers.)
\begin{gather*}
c(\Gamma)\leq m+1\\
\xi + E\leq m+1\\
\xi+m\leq 2m+1-E =\chi(\Gamma )\\
\xi+m\leq b_0 - b_1\\
b_0-\xi\geq m+b_1
\end{gather*}
Hence there are at least $m$ connected components of $\Gamma$ which have
no framed linked vertices. Picking a vertex from each such component
yields the desired $m$ free vertices.

In order to proceed, we need the following
interesting lemma. Let
$\mathcal V$ be the set of vertices of our grope.

\begin{lemma} \label{puppy}
Suppose $F\subset {\mathcal V}$ is a set of $m$ free vertices, $F = \{
v_1, \ldots , v_m\}$. We can assume $c( \Gamma\backslash {\rm star}F)
=0$. That is, if we remove $F$ and all edges which hit $F$ from $\Gamma$,
the complexity of the resulting graph is $0$.
\end{lemma}
\emph{[Proof]}

Suppose otherwise. Let $G$ be a class $2m+1$ grope with a set of $m$ free
vertices, $F$, contradicting the claim, with $c(\Gamma \backslash {\rm
star} F)$ minimal. By hypothesis this number is bigger than zero. Let $S
= \{ s_1,\ldots ,s_{m+2} \}$ be the scheme in which $s_1, \ldots,
s_{m-1}$ are type II moves trivializing the $v_1,\ldots ,v_{m-1}$ handles
supported between separating planes. $s_m, s_{m+1}$ are the in and out
move respectively on the $v_m$ handles. These two moves are supported in a
neighborhood of the $v_m$ handles with caps, which is separated from
the $v_1,\ldots ,v_{m-1}$ handles by hyperplanes, and so is disjointly
supported from the type II moves. Finally, $s_{m+2}$
is a type I move which reduces $c(\Gamma \backslash {\rm star} F)$. It
is possible that ${\rm supp} (s_{m+2})\cong \amalg D^3$ is not disjoint
from the other moves, since the separating planes may intersect this
disjoint union of balls. However, since $s_{m+2}$ is only reducing
complexity away from
$v_1,\ldots ,v_m$, at least the handles $v_1,\ldots , v_m$ do not hit
${\rm supp} (s_{m+2})$. But then the separating planes are easily pushed
out of
${\rm supp} (s_{m+2})$ using the balls to guide the isotopy, say. It is
then an easy matter to separate these balls from the other moves.

So $\sum_{\sigma \subset S} (-1)^{|\sigma |}[\partial G]_{m+1} = 0$,
and let us see what this says. In preparation, let us suppose that $G$
is formed by attaching the gropes $H^\prime$ and $H^{\prime\prime}$ to
the dual bands of the bottom stage, thereby partitioning $\mathcal V$
into two nonempty sets ${\mathcal V}_{H^\prime}$ and ${\mathcal
V}_{H^{\prime\prime}}$. Suppose without loss that $v_m\in {\mathcal
V}_{H^{\prime}}$. Let $S_{H^\prime}$ and $S_{H^{\prime\prime}}$ partition
$\{ s_1,\ldots , s_{m-1} \}$ into two sets in the obvious way. Let $S_I =
\{ s_m,s_{m+1} \}$ and $S_C = \{ s_{m+2}\}$. 

Note that we can assume
$s_{m+2}$ reduces $c(\Gamma \backslash {\mathcal V}_{H^\prime})$ since if
this complexity were zero, then ${\mathcal V}_{H^{\prime\prime}}$ would
have no edges hitting it, (and no framed linked vertices). By the
earlier stated assumption that the height function separates the two
halves of the grope $H^{\prime}$ and $H^{\prime\prime}$, the handles on
the $H^{\prime\prime}$ half all bound disks, implying of course that the
grope is trivial contradicting that $G$ is a counterexample. Thus we can
assume some complexity not contained wholly within the $H^\prime$ half,
and without loss
$s_{m+2}$ reduces this. 

We are now in a position to describe what happens under the various
combinations of moves from $S_{H^\prime}, S_{H^{\prime\prime}}, S_I$ and
$S_C$, with the initial assumption that neither $S_{H^\prime}$ nor
$S_{H^{\prime\prime}}$ is empty. In the following list of cases, case
$i$ refers to a set of moves, $\sigma$, which hits $i$ of the above $4$
sets.

\hspace{-1em}\emph{Case 0}

This is the empty move yielding {\boldmath $\partial G$}.

\hspace{-.5em}\emph{Case 1}

By our previous analysis of the handle trivializing moves, if $\sigma
\subset S_{H^\prime}$ or $\sigma \subset S_{H^{\prime\prime}}$, $\partial
G_\sigma$ is the {\bf unknot}. $K_{s_{m+2}}$ has less of the approriate
complexity so by minimality $[K_{s_{m+2}}]_{m+1}={\bf 0}$. The left over
terms are the ones gotten from the in/out trick: doing both of $s_m,
s_{m+1}$ is the {\bf unknot}, while $K_{s_m},K_{s_{m+1}}$ are
{\boldmath$ \partial H^\prime$ and
$ \rho(\partial H^\prime )$}.

\hspace{-1em}\emph{Case 2}

\hspace{-.5em}$\sigma$ hits $S_{H^\prime},S_{H^{\prime\prime}}$ :
{\bf unknot}.

\hspace{-.5em}$\sigma$ hits $S_{H^\prime},S_I$: $S_{H^\prime}$ trivializes
some handles, and then $s_m$ or $s_{m+1}$ give $H^{\prime}$ with
trivialized handles, an {\bf unknot}. Doing both the in and out move also
yields an {\bf unknot}.

\hspace{-.5em}$\sigma$ hits $S_{H^\prime},S_C$: $S_{H^\prime}$ trivializes
handles of the grope $G_{s_{m+2}}$ yielding an {\bf unknot}.

\hspace{-.5em}$\sigma$ hits $S_{H^{\prime\prime}}, S_I$:
$S_{H^{\prime\prime}}$ gives some grope with the $H^\prime$ half
unaltered. Doing one move from $S_I$ then gives the $H^\prime$ half.
Specifically,
 {\boldmath$
\sum_{\emptyset\neq\tau\subset S_{H^{\prime\prime}}} (-1)^{|\tau |}\{
[\partial H^\prime ]_{m+1}+[\rho(\partial H^\prime )]_{m+1}\} $}. Again
if we do both $s_m$ and $s_{m+1}$ the result is obviously an {\bf unknot}.

\hspace{-.5em}$\sigma$ hits $S_{H^{\prime\prime}}, S_C$: {\bf unknot}.

\hspace{-.5em}$\sigma$ hits $S_I,S_C$: $S_C$ gives some grope with the
$H^\prime$ half unaffected. So as in a previous case we get
 {\boldmath $
\partial H^\prime + \rho (\partial H^{\prime\prime})$} .

\hspace{-1em}\emph{Case 3}

\hspace{-.5em}$\sigma$ hits $S_{H^{\prime\prime}}, S_I, S_C$:
$S_{H^{\prime\prime}}, S_C$ give a grope with $H^\prime$ half intact,
and so as in two of the previous cases we get, adjusting the sign to
include the $s_{m+2}$ move, {\boldmath $$ \sum_{\emptyset \neq
\tau\subset S_{H^{\prime\prime}}} (-1)^{|\tau |+1} \{ [\partial H^\prime
]_{m+1}+[\rho(\partial H^\prime )]_{m+1}\}.$$ }

\hspace{-.5em}$\sigma$ hits $S_{H^\prime}, S_I, S_C$: {\bf unknot}.

\hspace{-.5em}$\sigma$ hits $S_{H^\prime},S_{H^{\prime\prime}}, S_C$:
{\bf unknot}.

\hspace{-.5em}$\sigma$ hits $S_{H^\prime}, S_{H^{\prime\prime}}, S_I$:
{\bf unknot}.

\hspace{-1em}\emph{Case 4}

This involves doing at least one move from each group and is an
{\bf unknot}. 
\vspace{1em}

We conclude 
\begin{gather*}
\sum_{\sigma\subset S} (-1)^{|\sigma |} [\partial
G_\sigma]_{m+1}= [\partial G] - [\partial H^\prime ]- [\rho (\partial
H^\prime )] +
\sum_{\emptyset\neq \tau\subset S_{H^{\prime\prime}}}(-1)^{|\tau |}\{
[\partial H^\prime] + [\rho (\partial H^\prime)]\}\\
 + [\partial H^\prime
] + [\rho (\partial H^\prime )] + \sum_{\emptyset\neq \tau\subset 
S_{H^{\prime\prime}}}(-1)^{|\tau +1
|}\{ [\partial H^\prime] + [\rho (\partial H^\prime)]\} = [\partial
G]_{m+1} = 0
\end{gather*} This is a contradiction. 

If $S_{H^\prime}=\emptyset$, then only cases leading to an
$m+1$-trivial knot are eliminated so the calculation still goes through.

If $S_{H^{\prime\prime}} = \emptyset$, then two nontrivial cases are
eliminated: the $S_{H^{\prime\prime}},S_I$ subcase of case 2 and the
$S_{H^{\prime\prime}}, S_I, S_C$ subcase of case 3. The calculation is
now $\sum_{\sigma\subset S} (-1)^{|\sigma |} [\partial G_\sigma ]_{m+1} =
[\partial G] - [\partial H^\prime ]-[\rho (\partial H^\prime )] +
[\partial H^\prime ] + [\rho (\partial H^\prime )] = 0$ which still
achieves the desired result $[\partial G]_{m+1} = 0$. 
$\Box$

Continuing the proof of theorem (\ref{mainthm}), recall we had found a
free set of $m$ vertices $F$. But the preceding lemma proves that
${\mathcal V}\backslash F$ can also be assumed free, this time of
cardinality $m+1$. Indeed we may assume that for any free $F^\prime$ of
cardinality $m$, ${\mathcal V}\backslash F^\prime$ is also free. This
actually implies $c(\Gamma )=0$ and therefore that $G$ is trivial, and we
are done: since $F$, ${\mathcal V}\backslash F$ are free, all
framed linked vertices have been eliminated. Suppose
${\mathcal V}\backslash F =
\{ w_1,
\ldots, w_{m+1}\}$. Let $F^\prime = \{ w_1,\ldots, w_m\}$. Then $F\cup \{
w_{m+1}
\}$ must be free, implying $w_m$ shares an edge with  no vertex in $F$.
Since it shared none with ${\mathcal V}\backslash F$, $w_{m+1}$ is in fact
isolated. But then by symmetry all of ${\mathcal V}\backslash F$ is
isolated. Since the only edges were between ${\mathcal V}\backslash F$
and $F$, there are no edges whatsoever. $\Box$

\newpage

\section{Showing the Bound is Sharp} \label{sharp}
In the following section we show that for all $n\geq 2$ there are knots
bounding gropes of class $n$ which are not $\lceil
\frac{n}{2}\rceil+1$-trivial. In fact, we find $K$ such that $J^{(\lceil
\frac{n}{2}\rceil+1)}_K(1)\neq 0$, where $J_K(t)$ is the Jones
polynomial. It is well known that the
$j$th derivatives of the Jones polynomial evaluated at $1$ are type $j$
invariants. Note that $J^{(m)}_\bullet(1)$ is not additive under connect
sum (primitive), but is easily seen to be additive on $m-1$-trivial knots.

For this section, it is convenient to use a different graph than the one
we used previously.

\begin{definition}
Let $G$ be a grope of class $n$ in standard position with framed unlinked
handles bounding fixed caps. We define the graph $\Gamma_\Delta (G)$ as
follows. The vertices as before correspond to the $V_i$, the
$n$ collections of handles satisfying the trivialization property.
We put in an edge between $V_i$ and $V_j$ if the corresponding caps
intersect.
\end{definition}

Note that type II moves on a free set of vertices have their obvious
analog in this setting: we make the moves by using the cap to guide the
homotopy. We call these \emph{type $II_\Delta$} moves, for clarity. The
moves are then obviously disjointly supported since the caps are
hypothesized to be disjoint.

We prove the following statement inductively:

\begin{theorem}\label{exthm}
For all even $n$, there is a grope $G$ of class $n$ with all the cores
$V_i$ unknotted, such that the corresponding graph $\Gamma_\Delta$ has no
first homology, and such that each vertex has valence less than or equal
to 2. Further the edges ending in valence 1 vertices correspond to finger
moves in the following sense: Let $\{ V_i,V_j\}$ be an edge with $V_i$
valence $1$. Suppose $V_i = \{ v_i^1,\ldots , v_i^{l_i}\}, V_j= \{ v_j^1,
\ldots ,v_j^{l_j}\}$. Then we insist the $\overline{v}_i^k$ (resp.
$\overline{v}_j^k$) bound disks $\Delta_i^k$ (resp. $\Delta_j^k$)
such that each $\Delta_i^k$ inersects exactly one $\Delta_j^k$ in a
single clasp singularity. Each $\Delta^k_j$ disk is hit at most once by
the $\Delta_i^k$ disks.

 This grope satisfies 
$J_K^{(\lceil
\frac{n}{2}\rceil+1)}(1)\neq 0$.
\end{theorem}

This is sufficient for our purposes since it also implies the odd case.
Just think of a grope of class $2m$ as a grope of class $2m-1$ by
ignoring one of the top stages. Since $\lceil\frac{2m}{2}\rceil =
\lceil \frac{2m-1}{2}\rceil = m$, any example of class $n=2m$ with
$J^{\lceil \frac{n}{2}\rceil+1}(1)\neq 0$ is also an example as a
class $n=2m-1$ grope. 

\emph{[ Proof ]}

If $n>2$, note that a graph satisfying the induction hypothesis will have
$\lceil
\frac{n}{2}\rceil$ free vertices and
 two
 special disjoint edges each containing a vertex which is not
contained in any other edge. To see this note that such a graph is
contained in a graph which is homeomorphic to an interval, the free
vertices being alternating vertices in this graph, and the
special edges being the edges at the ends of the interval. 

\begin{figure}
\begin{center}
\epsfig{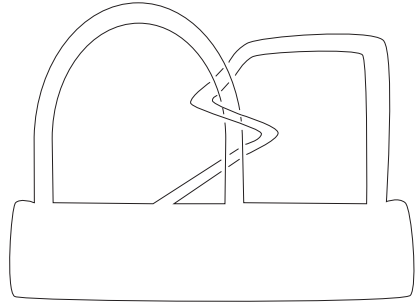}
\end{center}
\caption{This knot has $J(t)=2-t+t^2-2t^3+t^4-t^5+t^6$, as calculated by
Knotscape, with
$J^{(2)}(1)=12\neq 0$. \label{basecase1}}
\end{figure}

We need the base cases $n=2$, see figure (\ref{basecase1}), and $n=4$, the
second of which we defer to the end since we build the $n=4$ example from
the
$n=2$ example using a construction of the proof. The $n=2$ example does
not suffice because in order to get the induction going we need the graph
to have at least 2 edges. 

\begin{figure}
\begin{center}
\epsfig{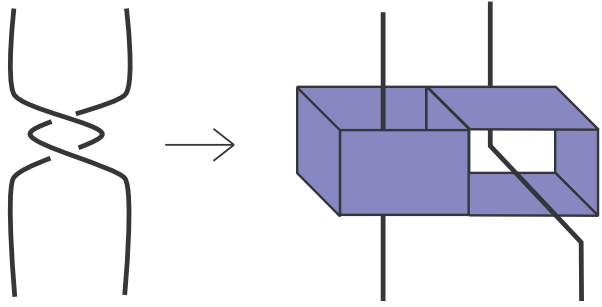}
\end{center}
\caption{Unlinking two handles and relinking them with the torus
$T_\alpha$.\label{figb}}
\end{figure}

Now assume $G$ is such a grope satisfying the statement of theorem
\ref{exthm} for
$n=2m$. Suppose the two special edges have endpoints $V_i, V_j$ and
$V_l,V_m$ respectively, with
$V_i$ and $V_l$ the ``dangling'' vertices. Take the edge $\{ V_i,V_j\}$ in
$\Gamma (G)$ and delete it, that is unlink the corresponding pairs
of handles of
$G$. Link each pair of these handles with a punctured torus $T_\alpha$ as
in figure \ref{figb}.

Notice that when the pairs of handles are pushed across each other to
relink, the boundary of the punctured tori will bound a symmetric surgery
disk. Denote by $\tilde{G}$ the grope $G$ modified as in figure
\ref{figb}, that is with the edge $\{ V_i,V_j\}$ deleted. Connect the 
boundaries of the punctured tori,
$T_\alpha$, with the bottom stage of $\tilde{G}$ by some bands disjoint
from the rest of the
$T_\alpha$ and from $\tilde{G}$ and also disjoint from the various caps
associated to all the vertices of
$\tilde{G}$. Call this new grope
$H_{ij}$. If $J^{(m+1)}_{H_{ij}}\neq 0$ let $H=H_{ij}$ and proceed.

\begin{figure}
\begin{center}
\epsfig{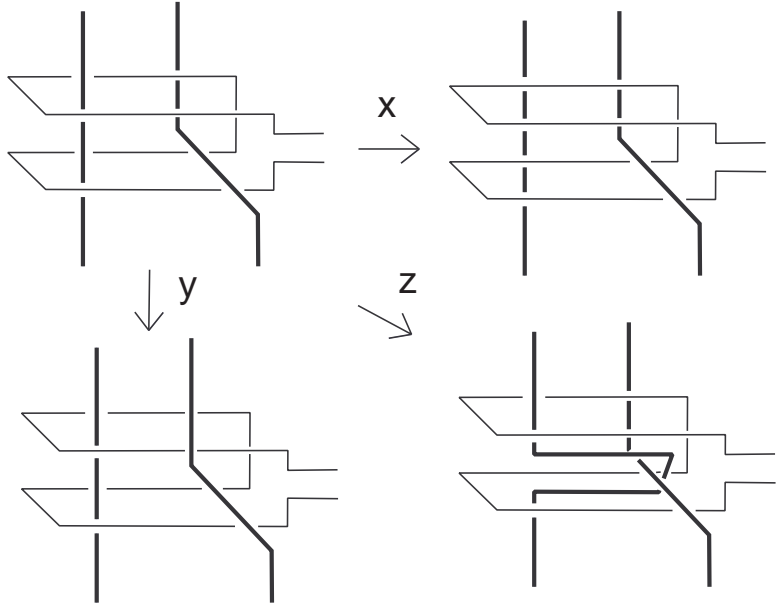}
\end{center}
\caption{The $x$, $y$ and $z$ moves.}
\label{figc}
\end{figure}

Otherwise, carry out the same procedure for the edge $\{ V_{l}, V_m\}$.
If this also fails, i.e. $J^{(m+1)}_{H_{lm}}(1)= 0$, we
form the grope $H_{lmij}$, which is the grope gotten from doing
the above procedure to \emph{both} edges. Consider the scheme $\{
s_1,\ldots ,s_{m-2}, x_{ij}, y_{lm} ,z_{ij},z_{lm} \}$, where the $s_i$
 are type II moves trivializing the $i$ handles corresponding to vertices
in the complement of the special edges, and where the $x$'s, $y$'s and
$z$'s are given on the corresponding $T_\alpha$ as pictured in
figure (\ref{figc}). The added torus, $T_\alpha$, has two bands
$x_\alpha$ and $y_\alpha$ each linking a handle of $\tilde{G}$ exactly
once. The move $x$ has the effect of removing the linkage of the
appropriate handle with $x_\alpha$ for all $\alpha$, whereas $y$ has the
corresponding effect on the $y_\alpha$. Indeed,  as the reader may verify,
doing any combination of these three moves $x,y,z$ on a particular
$T_\alpha$  causes this added torus to become compressible.
There is a choice in which bands are called $x_\alpha$ and which
$y_\alpha$. Denote by $x_\alpha$ those bands which links $V_i$ or $V_m$,
and the $y_a$ are then those bands linking $V_j$ or $V_l$.

Now
\begin{gather}
-[H_{lmij}]_{m+1}=\sum_{\emptyset\neq\sigma\subset\{
x_{ij},y_{lm},z_{ij},z_{lm}\}} (-1)^{|\sigma |} [H_{lmij\sigma}]_{m+1}
\label{har}\\ =\sum_{\emptyset\neq\sigma\subset\{ z_{ij},z_{lm}\}}
[H_{lmij}]_{m+1}\label{harhar}\\
=-[H_{lm}]_{m+1}-[H_{ij}]_{m+1}+[G]_{m+1} \label{harharhar}
\end{gather}
Here (\ref{har}) follows since doing any of the $s_i$ even in conjunction
with other moves in the scheme will cause the $\tilde{G}$ half to
trivialize, followed by the $t_\alpha$. (\ref{harhar}) follows since
doing either of $x_{ij}$ or $y_{lm}$ causes there to be a trivial group
of handles corresponding to a vertex in the $\tilde{G}$ half which then
trivializes the grope. 
Finally, (\ref{harharhar}) follows since doing $z_{ij}$, say, relinks the
$ij$ handles while causing the appropriate $T_\alpha$ to compress,
leaving the $T_\alpha$ linking with the $l,m$ handles, i.e. $H_{lm}$.

From (\ref{harharhar}) we could immediately conclude (\ref{ibcrootbeer}),
despite the fact that $J^{(m+1)}_\cdot (1)$ is not in general additive in
view of section 1.2, $(\ref{altsum})$. However, since we need it later
anyway, we will prove that $G,H_{lmij}, H_{lm},$ and $H_{ij}$ are all
$m$-trivial. Well $G$ is $m$-trivial by the main theorem. $H_{ij}$ is
$m$-trivial: Let
$s_1,\ldots ,s_{m-1}$ be type $II_\Delta$ moves corresponding to free
vertices in the complement of $\{ V_i ,V_j\}$. (Their existence is proven
later.) Consider the scheme $S=\{ s_1,\ldots, s_{m-1}, x,y\}$. Obviously,
any subset of these trivializes $H_{ij}$. Symmetrically $H_{ij}$ is
$m$-trivial. But (\ref{harharhar}) indicates that $H_{ijlm}$ is $m+1$,
hence $m$, equivalent to a sum of $m$-trivial knots. It is therefore
$m$-trivial itself. 
 
Thus 
\begin{gather}
J^{(m+1)}_{H_{lmij}}(1) = J^{(m+1)}_{H_{lm}}(1) + J^{(m+1)}_{H_{ij}}(1) -
J^{(m+1)}_G(1) \label{ibcrootbeer}\\
=0+0-J^{(m+1)}_G(1)\neq 0
\end{gather}
 We may let $H=H_{lmij}$. 

Recall that the $T_\alpha$ are connected via bands to $\tilde{G}$. We had
a lot of choice in choosing these and may assume they are organized as
follows. The $T_\alpha$ are band connect summed together to form $T$,
which is then connected by a band with $\tilde{G}$ which it links
geometrically.

We form a class
$n+2$ grope $K$ from
$H$ by plumbing as follows:
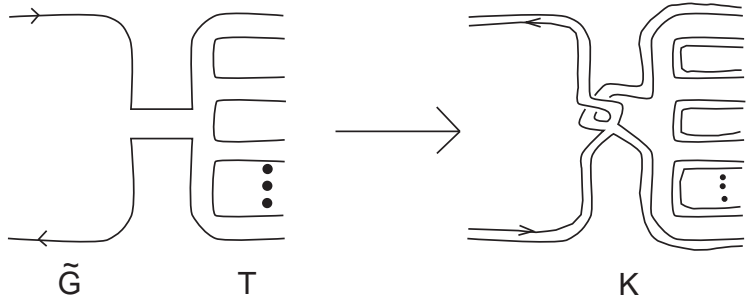
\begin{center}
\epsfig{file=fige ,%
        height=4cm}
\end{center}

\begin{figure}
\begin{center}
\epsfig{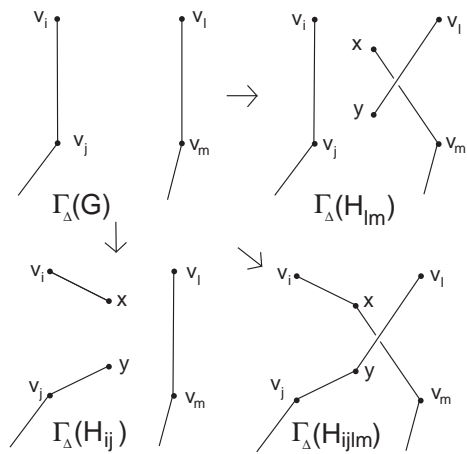}
\end{center}
\caption{The three possible moves in the construction of $H$.
}\label{graphmoves}
\end{figure}

That is, $K$ is formed by running a perpendicular annulus along
$\tilde{G}$ and one along $T$, and then plumbing these two annuli
together to get a punctured torus, the bottom stage of a new grope.
$K$ is a class $n+2$ grope, the bottom stage of which has a core
bounding a class $n$ grope which was gotten from $G$, and the dual core
of which bounds a class $2$ grope which is the connected sum of
the punctured tori, $T_\alpha$. We claim
$J^{(m+2)}_K(1)\neq 0$ which will complete the inductive step since
$\lceil\frac{n+2}{2}\rceil+1 =
\lceil\frac{n}{2}+1\rceil+1 = \lceil\frac{n}{2}\rceil+2 = m+2$. Let $x$
and $y$ denote half symplectic bases of $T$. Then the altered graph is as
in figure (\ref{graphmoves}). Suppose
$\{ s_i\}_{i=1}^{m}$ are $m$ free vertices on the $G$ half of $K$, none of
which is on an edge connected to $x$. Consider the scheme $S=\{
s_i\}_{i=1}^{m-1}\cup \{ {\rm in}x, {\rm out}x, {\rm in}s_m, {\rm out}
s_m\}$, where the $s_i$ are type II moves making the respective handles of
$K$ bound disks. As we know, if we do any of the $s_i$, then $K$
trivializes. Thus,
\begin{eqnarray}
-[K]_{m+2}=\sum_{\emptyset\neq \sigma \subset \{{\rm in}x_m,{\rm out}x_m
,{\rm in}x,{\rm out}x\}} (-1)^{|\sigma |}[K_\sigma ]_{m+2}  \\
=[H]+[\hat{H}]+[\rho (H)] + [\rho (\hat{H})] \label{5star}
\end{eqnarray}
Where \ref{5star} follows from lemma \ref{4term}. Let us
compute the relationship between $H$ and $\hat{H}$ in terms of the Jones
polynomial.
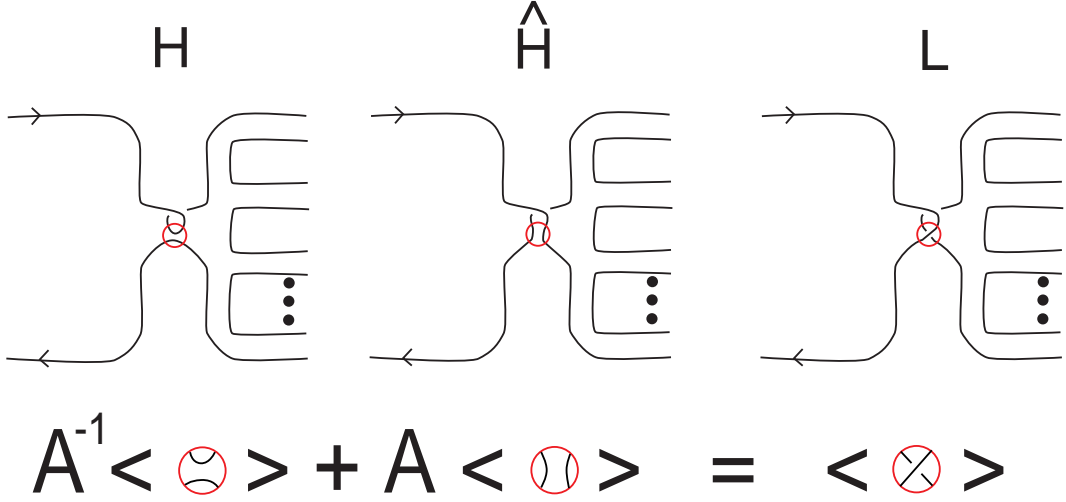
\begin{center}
\epsfig{file=figf ,%
        height=7cm}
\end{center}

That is $A<\hat{H}> + A^{-1}<H>=<L>$, where $<\bullet >$ denotes the
Kauffman bracket. Assume the writhe $w$ of the diagrams is zero away from
the pictured spots. Then $w(H)=1, w(\hat{H})=-1$ and $w(L)=0$. Thus $J_H=
(-A)^{-3} <H>$, $J_L = <L>$ and $J_{\hat{H}}=(-A)^3<H>$. This implies the
following relation, where we make the substitution
$A^{-2}=t^{\frac{1}{2}}$,
\begin{equation}
-t^{\frac{1}{2}}J_{\hat{H}}(t)-t^{-\frac{1}{2}}J_H(t)=J_L(t)
\end{equation}
Setting $u=t^{\frac{1}{2}}$, we get
\begin{eqnarray}
uJ_{\hat{H}}(u)-u^{-1}J_H(u)=J_L(u) \label{star}
\end{eqnarray}

Note that $\hat{H}$ is also $m$-trivial. We now analyze the triviality of
$L$.
\begin{claim}
$J^{(k)}_L(1)=J^{(k)}_{Unlink}(1)$ for all $k\leq m+2$.
\end{claim}

\emph{[Proof]}

First, choose $m$ free vertices $x_1,\ldots, x_m$ on the $\tilde{G}$ half
of $L$, none of which shares an edge with the $x$ vertex in
$\Gamma_\Delta (K)$. The possible ways that we altered the graph are
listed in figure
\ref{graphmoves}. Each is contained in a graph which spans all vertices
and is homeomorphic to a line as remarked earlier, which we can always
choose with $x$ as the second vertex from an endpoint. It is then obvious
we can choose $m$ vertices which do not share an edge with $x$.  

\begin{figure}
\begin{center}
\epsfig{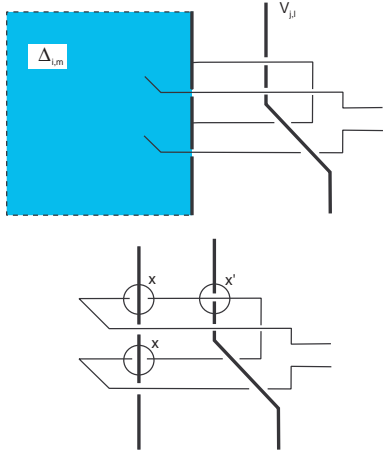}
\end{center}
\caption{The $x$ and $x^\prime$ moves.}
\label{listrivial}
\end{figure}

Without loss, the local pictures of the tori $T_\alpha$ look like figure
\ref{listrivial}, with the $V_i$ or $V_m$ cap $\Delta_i$ or $\Delta_m$
hitting $T_\alpha$ as indicated.

Consider the two indicated sets of crossing changes, where `$x$' is our
old friend. Doing $x'$ makes all $V_i$ bound their caps in the complement
of $T$. $\tilde{G}$ is thus isotopic to the unknot in the complement of
$T$ since the isotopy is supported in a neighborhood of $\tilde{G}\cup
V_i {\rm caps}$. We are then left with an unlink since $T$ is unknotted
when $\tilde{G}$ is removed. Similarly doing $x$ leaves an unlink of two
components. Consider the scheme $S=\{ s_1,\ldots ,s_{m-1}, {\rm in}x_m ,
{\rm out}x_m,x,x'\}$, where $s_i$ are type $II_\Delta$ moves trivializing
the
$x_i$. Doing such a type $II_\Delta$ move also yields an unlink, because
after the $x_i$ handles are trivialized, they bound caps in the
complement of
$T$ and so the previous argument goes through. Similarly, the in and out
moves on $x_m$ trivialize the knot in a neighborhood of $\tilde{G}\cup
x_m {\rm cap}$. Indeed doing any combination of moves in $S$ gives an
unlink. Since $|S|=m+3$, it follows that
$J^{(k)}_L(1)=J^{(k)}_{unlink}(1)$ for all $k\leq m+2$. (Recall
$J^{(k)}$ is a type $k$ link invariant. $\Box$

 Applying $(\frac{d}{dt})^{m+2} = 
(\frac{du}{dt})^{m+2}(\frac{d}{du})^{m+2}$ to both sides of (\ref{star}),
and using claim 1, we get 
\begin{eqnarray}
-\left(\frac{d}{du}\right)^{m+2} (uJ_{\hat{H}}(u))  (1) -
\left(\frac{d}{du}\right)^{m+2}(u^{-1}J_H(u)) (1) =
\left(\frac{d}{du}\right)^{m+2}(J_{{\rm unlink}}(u)) (1) \label{dagg}
\end{eqnarray}

To evaluate each of these, note $J_{\rm unlink} = -A^{-2}-A^2 = -u -
u^{-1}$. So the right hand side of (\ref{dagg}) is equal to
$-(\frac{d}{du})^{m+2}(u^{-1})|_1$. Also 
\begin{gather*}
\left(\frac{d}{du}\right)^{m+2}(u
J_{\hat{H}}(u))(1)=\sum_{k=0}^{m+2} \binom{m+2}{k}\left(\frac{d}{du}
\right)^k(u)|_1 \left(\frac{d}{du}\right)^{m+2-k}(J_{\hat{H}})_1 =\\
 J^{(m+2)}_{\hat{H}}(1) +
(m+2)J^{(m+1)}_{\hat{H}}(1).
\end{gather*}
 Finally,
\begin{gather*}
\left(\frac{d}{du}\right)^{m+2}(u^{-1}J_H) \\
=\sum_{k=0}^{m+2}\binom{m+2}{k}\left(\frac{d}{du}\right)^k(u^{-1})|_1
\left(\frac{d}{du}\right)^{m+2-k}(J_H)(1)\\
= J^{(m+2)}_H(1) +(m+2)\frac{d}{du}(u^{-1})(1)J^{(m+1)}_H(1) +
\left(\frac{d}{du}\right)^{m+2}(u^{-1})|_1 J^{(0)}_H(1)\\
= J^{(m+2)}_H(1) - (m+2)J^{(m+1)}_H(1) +
\left(\frac{d}{du}\right)^{m+2}(u^{-1})|_1
\end{gather*}
Thus equation (\ref{dagg}) becomes
\begin{eqnarray}
-J^{(m+2)}_{\hat{H}}(1) - (m+2)J^{(m+1)}_{\hat{H}}(1)    
-J^{(m+2)}_H(1) + (m+2)J^{(m+1)}_H(1) =0 \label{4star}
\end{eqnarray}

I claim that $J^{(m+2)}_{\hat{H}}(1) + J^{(m+2)}_H(1)\neq 0$. 
Otherwise, (\ref{4star}) implies that $J^{(m+1)}_{\hat{H}}(1)=
J^{(m+1)}_H(1)
$. Consider (\ref{5star}). It implies that $J^{(k)}_{K^{-1}}(1) =
J^{(k)}_{H\#\rho (H)\#\hat{H} \#\rho (\hat{H})}(1)$ for all $k\leq m+2$.
But
$K$ is $m+1$-trivial by the main theorem, hence $J^{(m+1)}_{H\#\rho
(H)\#\hat{H} \#\rho (\hat{H})}(1) = 0$. Since $H$ and $\hat{H}$ are $m$
trivial, and $J$ is invariant under $\rho$, this implies
$2J^{(m+1)}_{H}(1) + 2J^{(m+1)}_{\hat{H}}(1) = 0.$ This would imply
$J^{(m+1)}_H(1)=J^{(m+1)}_{
\hat{H}}(1)=0$, contradicting our choice of $H$.
 
So $J^{(m+2)}_{\hat{H}}(1) + J^{(m+2)}_H(1) \neq  0$. Note $J_{H\#\rho
(H) \# \hat{H} \# \rho (\hat{H} )}= J_H^2 J_{\hat{H}}^2$. So

\begin{gather}
J^{(m+2)}_{H \# \ldots \# \rho(\hat{H})}(t) =\\
\left(\frac{d}{dt}\right)^{m+2} \left(J_H^2(t)J_{\hat{H}}^2(t)\right)
=\sum_{k=0}^{m+2}\binom{m+2}{k}\left(\sum_{l=0}^{k}\binom{k}{l}
\left( \frac{d}{dt}\right)^l J_H \left(\frac{d}{dt}\right)^{k-l}J_H\right)
\times
\\ \left( \sum_{l=0}^{m+2-k}
\binom{m+2-k}{l}\left(\frac{d}{dt}\right)^lJ_{\hat{H}}
\left(\frac{d}{dt}\right)^{m+2-k-l}J_{\hat{H}}\right) . \label{yuppy}
\end{gather}
   
In order for $(\frac{d}{dt})^{k-l}J_H(1)\neq 0$, $k-l\in \{ 0,m+1,
m+2\}$, which means that either ($k=l$) or ($l=0$ and $k=m+1$) or ($k=m+2
$  and $ (l=0 $ or  $l=1))$. In order for
$(\frac{d}{dt})^{l}J_H(1)\neq 0$, $l\in \{ 0, m+1, m+2\}$. So the only
potentially nonzero terms that arise upon evaluating (\ref{yuppy}) occur
when 
$k=l=0,m+1,m+2$ or
$k=m+1,l=0$ or
$k=m+2, l=0$. So we get, for $k=0$,

\begin{gather*}
(J_H(1)J_H(1))\left( \sum_{l=0}^{m+2} \binom{m+2}{l}
\left(\frac{d}{dt}\right)^l\mid_1 J_{\hat{H}}(t)
\left(\frac{d}{dt}\right)^{m+2-l}\mid_1 J_{\hat{H}}(t)\right) =\\ 1\cdot 1
(J_{\hat{H}}(1) J^{(m+2)}_{\hat{H}}(1) + J^{(m+2)}_{\hat{H}}(1)
J_{\hat{H}}(1)) =\\ 2J^{(m+2)}_{\hat{H}}(1).
\end{gather*}

For $k=m+1$, we get
\begin{gather*}
(m+2)\left(J_H(1)J^{(m+1)}_H(1)\right)\left(\sum_{l=0}^{1} 
\binom{1}{l} \left(\frac{d}{dt}\right)^l
\mid_1 J_{\hat{H}}(t) \left(\frac{d}{dt}\right)^{1-l}\mid_1
J_{\hat{H}}(t)\right)\\ =(m+2)J^{(m+1)}_H(1)\cdot 0 = 0
\end{gather*}

For $k=m+2$, we get 
\begin{gather*}
\left( J_H(1)J^{(m+2)}_H(1) + J^{(m+2)}_H(1)J_H(1)\right)
\cdot (J^{(0)}_{\hat{H}}(1)J^{(0)}_{\hat{H}}(1))=\\
2J^{(m+2)}_H(1).
\end{gather*}

So, equation (\ref{5star}) implies
$2J^{(m+2)}_H(1)+2J^{(m+2)}_{\hat{H}}(1) = J^{(m+2)}_{K^{-1}}$. Thus,
since $J^{(m+2)}_H(1)+J^{(m+2)}_{\hat{H}}(1)
\neq 0$, we have shown that $J^{(m+2)}_{K^{-1}}(1)\neq 0$. In ${\mathcal
G}_{m+2}$, $K^{-1}+K = 0$, and since $K, K^{-1}$ are $m+1$-trivial
( $K^{-1}\# K\sim_{m+2} 0 \Rightarrow K^{-1}\#K\sim_{m+1}0$ and then
$K^{-1}\# K\sim_{m+1}K^{-1}\#0 \Rightarrow K^{-1}\sim_{m+1} 0$) , this
implies $J^{(m+2)}_{K^{-1}}(1)=-J^{(m+2)}_K(1)\neq 0$, establishing the
key property of the inductive statement. 

We must also show that the graph of $K$ has no cycles and has valences
less than or equal to 2, and that each edge corresponds to a finger move.
 When one forms $K$ from $H$ the handle
pattern is the same.

\begin{figure}
\begin{center}
\epsfig{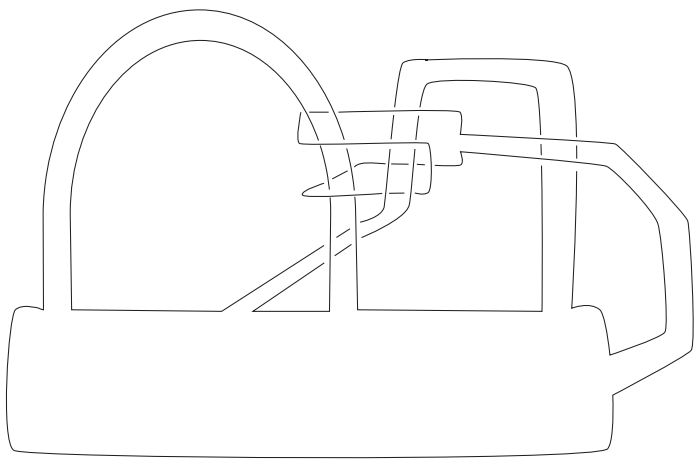}
\end{center}
\caption{In order for the induction to get started $J^{(2)}(1)$ should
not be zero on the above knot.}\label{dim4ctreg}
\end{figure}

To see that the two graph properties are the same, notice $K$ was
constructed with one of the two following moves on the graph. The first
one is to take a special edge, delete it, and then connect two added
vertices to the endpoints of the deleted edge. The second is to delete
two edges, add two vertices, and add four edges between the added
vertices.  The change is as diagrammed in figure \ref{graphmoves} and
preserves the properties we want.
The fact that the edges correspond to finger moves follows since in the
added surface $T$, each added torus has a handle which links a handle of
$G$ once in a single clasp.
Finally we must exhibit an example for the case $n=4$. The only possible
problem with the above induction would be if when forming $H$ as pictured
in figure \ref{dim4ctreg}, $J^{(2)}_H(1)=0$. However, using Knotscape to
calculate the Jones Polynomial, we see
$J_H(t)=t^{-4}-2t^{-3}+3t^{-2}-4t^{-1}+5-4t+3t^2-2t^3+t^4$, and one may
calculate $J^{(2)}_H(1)=12\neq 0$.

$\Box$
\newpage
\footnotesize
{\bf Bibliography}

[BL]  J. Birman and X.S. Lin, \emph{Knot polynomials and
Vassiliev's invariants}, Invent. Math. {\bf 111} (1993), 225-270

\vspace{.5 em}
[FT] M. Freedman and P. Teichner, \emph{ 4 Manifold Topology II: Dwyer's
Filtration and atomic surgery kernels}, Invent. Math. (1995)

\vspace{.5em}
[G] M. Gusarov, \emph{On $n$-equivalence of knots and invariants of
finite degree}, Topology of manifolds and varieties, Adv. Soviet Math.,
vol 18, Amer. Math. Soc., 1994, 173-192

\vspace{.5em}
[LK] X.S. Lin and E. Kalfagianni, \emph{Regular Seifert Surfaces and
Vassiliev Knot Invariants}, preprint available on X.S. Lin's homepage

\vspace{.5em}
[NS] K.Y. Ng and T. Stanford, \emph{On Gusarov's groups of knots},
Math.Proc. Cambridge Philos. Soc., to appear

\vspace{.5em}
[V] V.A. Vassiliev, \emph{Complements of discriminants of smooth maps:
topology and applications}, Trans. of Math. Mono. {\bf 98}, Amer. Math.
Soc., Providence, 1992

\end{document}